# Long time behavior of solutions to the 2D Boussinesq equations with zero diffusivity

Igor Kukavica
Weinan Wang

Monday 11$^{\text{th}}$ November, 2019

Department of Mathematics
University of Southern California
Los Angeles, CA 90089
e-mails: kukavica@usc.edu, wangwein@usc.edu

## Abstract

We address long time behavior of solutions to the 2D Boussinesq equations with zero diffusivity in the cases of the torus, $\mathbb{R}^2$, and on a bounded domain with Lions or Dirichlet boundary conditions. In all the cases, we obtain bounds on the long time behavior for the norms of the velocity and the vorticity. In particular, we obtain that the norm $\|(u, \rho)\|_{H^2 \times H^1}$ is bounded by a single exponential, improving earlier bounds.

## 1 Introduction

We consider the asymptotic behavior of solutions to the Boussinesq equations without diffusivity

$$u_t - \Delta u + u \cdot \nabla u + \nabla \pi = \rho e_2 \tag{1.1}$$

$$\rho_t + u \cdot \nabla \rho = 0 \tag{1.2}$$

$$\nabla \cdot u = 0 \tag{1.3}$$

in a bounded domain $\Omega \subseteq \mathbb{R}^2$, $\mathbb{T}^2$, and $\mathbb{R}^2$. Here, $u$ is the velocity satisfying the 2D Navier-Stokes equations [CF, DG, FMT, R, T1, T2, T3] driven by $\rho$, which represents the density or temperature of the fluid, depending on the physical context. Also, $e_2 = (0, 1)$ is the unit vector in the vertical direction.

Recently, there has been a lot of progress made on the existence, uniqueness, and persistence of regularity, mostly in the case of positive viscosity and vanishing diffusivity, considered here, while the same question with both vanishing viscosity and diffusivity is an important open problem. The initial results on the global existence in the regularity class have been obtained by Hou and Li [HL], who proved the global existence and persistence in the class $H^s \times H^{s-1}$ for integer $s \geq 3$. Independently, Chae [C] considered the class $H^s \times H^s$ and proved the global persistence in $H^3 \times H^3$. The class $H^s \times H^{s-1}$ has subsequently been studied in the case of a bounded domain, where Larios et al proved in [LLT] the global existence and uniqueness for $s = 1$ and then by Hu et al, who proved in [HKZ1] the persistence for $s=2$. The remaining range $1 < s < 3$ was then resolved in [HKZ2] in the case of periodic boundary conditions. For other works on the global existence and persistence in Sobolev and Besov classes, see [ACW, BS, BrS, CD, CG, CN, CW, DP1, DP2, DWZZ, HK1, HK2, HS, KTW, KWZ, LPZ].

In a recent paper [J], Ju addressed the important question of long time behavior of solutions. He proved that in the case of Dirichlet boundary conditions on a bounded domain $\Omega$, the $H^2(\Omega) \times H^1(\Omega)$ norm grows at most as $Ce^{Ct^2}$, where $C > 0$ is constant. In the present paper, we consider this question for this and other boundary conditions. When the domain is finite, we prove that actually the $H^2 \times H^1$ norm is increasing as a single exponential. We conjecture that this bound is sharp. This is because it is not expected that the solutions of the Boussinesq equation decay. However, note that the rate of increase of the gradient of the density is bounded by the exponential integral of the $L^\infty$ norm of the gradient, i.e.,

$$\|\nabla \rho(t)\|_{L^2} \lesssim \exp\left(\int_0^t \|\nabla u(s)\|_{L^\infty}\,ds\right) \|\nabla \rho_0\|_{L^2},$$

cf. (2.37) below, and if $u$ is not decaying, we should expect the integral to be bounded from below by a constant multiple of $t$. In addition to the behavior of $\|(u,\rho)\|_{H^2 \times H^1}$, we also address the long time behavior of the vorticity. In the case of the torus, we find constant upper bounds for the vorticity and the gradient of the vorticity for all $L^p$ norms. This result relies on the uniform upper bound for $\|u\|_{H^2}$ established in [J] as well as on a Nash-Moser type result on the growth of the vorticity, stated as Lemma 2.2 below and which we believe is of independent interest.

The paper is structured as follows. In Section 2, we first address the case of periodic boundary conditions. In this case, the exponential bound for the gradient of the density is obtained by establishing a constant upper bound for $\|\nabla u\|_{L^p}$. For this purpose, we first obtain a uniform upper bound for all the $L^p$ norms of the vorticity, a result based on a Nash-Moser type iteration. To do the same for the gradient of the vorticity, it is not suitable to proceed with direct estimates. Instead, we recall the concept of the generalized vorticity $\zeta$ (cf. (2.20) below), which reduces the number of the derivatives in the density by one.

In Section 3, we consider the case of the unbounded domain $\mathbb{R}^2$. Here, the energy does not decay and in fact, the quantity $\|u(\cdot,t)\|_{L^2}$ grows linearly in time. Applying a similar procedure as in Section 2, we obtain $\|u(\cdot,t)\|_{H^2} = \mathcal{O}(t^{1/2})$ as well as an information on the growth of $\|\rho\|_{H^1}$. In addition, we obtain upper bounds for $\|\omega\|_{L^p}$ and $p^{-3/2}\|\nabla \omega\|_{L^p}$, which are uniform in $p$.

In the final two sections, we address the case of a smooth bounded domain with either Lions or Dirichlet boundary conditions. For the Lions boundary conditions, we obtain $\|\nabla \rho\|_{L^2} \leq Ce^{Ct}$, using a different technique than the one for periodic boundary conditions. In addition, we obtain a uniform constant upper bound for $\|\omega\|_{L^p}$. Similarly, the last section contains the results in the case of Dirichlet boundary conditions, where we obtain an exponential upper bound for $\|\nabla \rho\|_{L^2} \lesssim e^{Ct}$, improving the main result in [J].

## 2 Long time behavior for periodic boundary conditions

In this section, we consider the Boussinesq system (1.1)–(1.3) in the case of the torus $\mathbb{T}^2$, i.e., assuming that $u$ and $\rho$ are 1-periodic. We assume for simplicity that $\int_{\mathbb{T}^2} u(\cdot,t) = 0$ for all $t \geq 0$; the general case can be addressed with the same methods; cf. Remark 2.3 below. The system is supplemented with the initial condition

$$(u(\cdot,0), \rho(\cdot,0)) = (u_0, \rho_0) \in H^2(\mathbb{T}^2) \times H^1(\mathbb{T}^2)$$

with $u_0$ divergence-free. By [HKZ1], there exists a global solution $(u(t), \rho(t))$ which belongs to $H^2 \times H^1$. Also, by [J], we have

$$\|u(t)\|_{H^2} \leq C, \qquad t \geq 0. \tag{2.1}$$

In the following statement, we provide an upper bound for the growth of the $\rho$ component of the norm $\|(u, \rho)\|_{H^2 \times H^1}$. Also, we establish a uniform upper bound on the quantities $\|\omega(\cdot, t)\|_{L^p}$ and $p^{-3/2}\|\nabla\omega(\cdot, t)\|_{L^p}$ for all $p \geq 2$.

**Theorem 2.1.** *Assume that $(u_0, \rho_0) \in H^2(\mathbb{T}^2) \times H^1(\mathbb{T}^2)$ satisfies $\nabla \cdot u_0 = 0$ and $\int_{\mathbb{T}^2} u_0 = 0$. Then we have*

$$\|\rho(t)\|_{H^1} \leq Ce^{Ct}, \qquad t \geq 0$$

*for a constant $C = C(\|u_0\|_{H^2}, \|\rho_0\|_{H^1})$. Moreover,*

$$\|\omega(t)\|_{L^p} \leq C, \qquad t \geq t_0, \qquad p \in [2, \infty]$$

*and*

$$\|\nabla\omega(t)\|_{L^p} \leq Cp^{3/2}, \qquad t \geq t_0, \qquad p \in [2, \infty), \tag{2.2}$$

*where $t_0 \geq 0$ depends on $\|u_0\|_{L^2}$.*

Note that (2.1) and (2.2) imply

$$\|u\|_{W^{2,p}} \leq Cp^{5/2}, \qquad t \geq t_0, \qquad p \in [2, \infty).$$

In the proof, we need the following statement on the long time behavior of solutions to the Navier-Stokes equations, which is of independent interest.

**Lemma 2.2.** *Consider the Navier-Stokes system*

$$u_t - \Delta u + u \cdot \nabla u + \nabla \pi = f$$
$$\nabla \cdot u = 0,$$

*supplemented with a divergence-free initial condition $u(\cdot, 0) = u_0 \in L^2(\mathbb{T}^2)$ such that $\int_{\mathbb{T}^2} u_0 = 0$ and $\int_{\mathbb{T}^2} f(\cdot, t) = 0$ for $t \geq 0$. If, for some $\lambda \geq 0$, we have*

$$\|f\|_{L^\infty([0,\infty), L^p(\mathbb{T}^2))} \leq p^\lambda M, \qquad 2 \leq p < \infty, \tag{2.3}$$

*where $M \geq 1$, then there exists $t_0 > 0$ depending only on $\|u_0\|_{L^2}$ such that*

$$\|\omega(\cdot, t)\|_{L^p} \leq CM, \qquad t \geq t_0, \qquad 2 \leq p \leq \infty, \tag{2.4}$$

*where $C$ is a universal constant. Moreover, for every $t_0 > 0$, there exists a constant $C$ depending only on $\|u_0\|_{L^2}$ and $t_0$ such that (2.4) holds.*

The proof uses ideas from [K, Lemma 3.1], where $\lambda = 0$ was considered. Lemma 2.2 is needed below with $\lambda = 1/2$.

*Proof of Lemma 2.2.* First, we prove (2.4) for some $t_0 > 0$, leaving the last assertion to the end of the proof. Without loss of generality, $M \geq 2$. The energy inequality reads

$$\frac{1}{2}\frac{d}{dt}\|u\|_{L^2}^2 + \|\nabla u\|_{L^2}^2 \leq \|f\|_{L^2}\|u\|_{L^2}, \tag{2.5}$$

from where, using the Poincaré inequality,

$$\frac{d}{dt}\|u\|_{L^2} + \frac{1}{C}\|u\|_{L^2} \leq \|f\|_{L^2}.$$

Applying the Gronwall inequality and shifting time, we may assume, without loss of generality, that

$$\|u(t)\|_{L^2} \leq CM, \qquad t \geq 0. \tag{2.6}$$

Note that the size of the time shift depends only on $\|u_0\|_{L^2}$ and $M$. Next, the vorticity $\omega = \nabla \times u$ satisfies

$$\omega_t - \Delta\omega + u \cdot \nabla\omega = \nabla \cdot F, \tag{2.7}$$

where $F = (F_1, F_2) = (f_2, -f_1)$. For $p = 2, 4, 8, \ldots$, define

$$\phi_p = \int \omega^p,$$

where all the integrals in this section are assumed to be over $\mathbb{T}^2$. First, the enstrophy inequality reads

$$\frac{1}{2}\phi_2' + \|\nabla\omega\|_{L^2}^2 \leq \|F\|_{L^2}\|\nabla\omega\|_{L^2} = \|f\|_{L^2}\|\nabla\omega\|_{L^2} \leq \frac{1}{2}\|f\|_{L^2}^2 + \frac{1}{2}\|\nabla\omega\|_{L^2}^2,$$

from where, using

$$\|\nabla\omega\|_{L^2}^2 \geq \frac{\|\omega\|_{L^2}^4}{\|u\|_{L^2}^2} = \frac{\phi_2^2}{\|u\|_{L^2}^2},$$

which follows from $\|\omega\|_{L^2} = \|\nabla u\|_{L^2} \leq \|u\|_{L^2}^{1/2}\|\Delta u\|_{L^2}^{1/2} = \|u\|_{L^2}^{1/2}\|\nabla\omega\|_{L^2}^{1/2}$, we obtain

$$\phi_2' + \frac{\phi_2^2}{C\|u\|_{L^2}^2} \leq \|f\|_{L^2}^2.$$

Therefore, by (2.6) and $\|f\|_{L^2} \lesssim M$,

$$\phi_2' + \frac{\phi_2^2}{CM^2} \leq CM^2,$$

and thus there exists a universal constant $t_1 \geq 0$ such that

$$\phi_2(t) \leq CM^2, \qquad t \geq t_1.$$

Now, let $p \in \{2, 4, 8, \ldots\}$. Testing the vorticity equation (2.7) with $\omega^{2p-1}$, we get

$$\frac{1}{2p}\phi_{2p}' + (2p-1)\int \omega^{2p-2}|\nabla\omega|^2 = \int \partial_j F_j \omega^{2p-1}$$

$$= -(2p-1)\int F_j \omega^{2p-2}\partial_j \omega \leq (2p-1)\|F\|_{L^{2p}}\|\omega^{p-1}\|_{L^{2p/(p-1)}}\|\omega^{p-1}\nabla\omega\|_2$$

$$\leq \frac{2p-1}{2}\int \omega^{2p-2}|\nabla\omega|^2 + Cp\|f\|_{L^{2p}}^2\|\omega^{p-1}\|_{L^{2p/(p-1)}}^2,$$

from where
$$\frac{1}{2p}\phi'_{2p} + \frac{2p-1}{2}\int \omega^{2p-2}|\nabla\omega|^2 \leq Cp\|f\|^2_{L^{2p}}\|\omega\|^{2p-2}_{L^{2p}}. \tag{2.8}$$

Using Nash's inequality, cf. [N, p. 936],
$$\|v\|_{L^2} \lesssim \|v\|^{1/2}_{L^1}\|\nabla v\|^{1/2}_{L^2} + \|v\|_{L^1} \tag{2.9}$$

with $v = \omega^p$, we get $\|\omega^p\|_{L^2} \lesssim \|\omega^p\|^{1/2}_{L^1}\|\nabla(\omega^p)\|^{1/2}_{L^2} + \|\omega^p\|_{L^1}$ whence $\|\omega\|^{4p}_{L^{2p}} \lesssim p^2\|\omega\|^{2p}_{L^p}\|\omega^{p-1}\nabla\omega\|^2_{L^2} + \|\omega\|^{4p}_{L^p}$. Therefore,
$$\|\omega^{p-1}\nabla\omega\|^2_{L^2} \geq \frac{\|\omega\|^{4p}_{L^{2p}} - C\|\omega\|^{4p}_{L^p}}{Cp^2\|\omega\|^{2p}_{L^p}}.$$

Applying this inequality on the second term in (2.8), we get
$$\frac{1}{2p}\phi'_{2p} + \frac{\phi^2_{2p} - C\phi^4_p}{Cp\phi^2_p} \leq Cp\|f\|^2_{L^{2p}}\phi^{(p-1)/p}_{2p}, \tag{2.10}$$

whence, by (2.3),
$$\phi'_{2p} + \frac{\phi^2_{2p} - C\phi^4_p}{C\phi^2_p} \leq Cp^{2+2\lambda}M^2\phi^{(p-1)/p}_{2p}.$$

Note that if
$$\phi_{2p} \geq C_0 \max\Big\{\phi^2_p, p^{2(1+\lambda)p/(p+1)}\phi^{2p/(p+1)}_p M^{2p/(p+1)}\Big\}, \tag{2.11}$$

then
$$\phi'_{2p} + \frac{\phi^2_{2p}}{C\phi^2_p} \leq 0,$$

which means that once $\phi_p$ is bounded, $\phi_{2p}$ is rapidly decreasing as long as it is sufficiently large. By increasing the constants, we may assume that
$$\phi_2(t) \leq C_0 M^2, \qquad t \geq t_1$$

and $C_0 \geq 1$. Denote $p_k = 2^k$, for $k \in \mathbb{N}$. Now, define recursively a sequence $M_1, M_2, M_3, \ldots$ such that
$$M_{k+1} = C_0 \max\Big\{p_k M^2_k, p^{2(1+\lambda)p_k/(p_k+1)}_k M^{2p_k/(p_k+1)}_k M^{2p_k/(p_k+1)}\Big\}, \qquad k = 1, 2, \ldots \tag{2.12}$$

(the reason for $p_k$ in front of $M^2_k$, comparing (2.12) with (2.11), is that it appears on the right side of (2.13) below). Also, let
$$M_1 = C_0 M^2.$$

We shall define a sequence $0 \leq t_1 \leq t_2 \leq \cdots$ such that
$$\phi_{2^k}(t) \leq M_k, \qquad t \geq t_k$$

with $\{t_k\}^\infty_{k=1}$ uniformly bounded. To construct this sequence, we proceed inductively, and assume that $t_k$ has been set. As long as $\phi_{2^{k+1}} \geq M_{k+1}$, we have
$$\phi'_{2^{k+1}} + \frac{\phi^2_{2^{k+1}}}{CM^2_k} \leq 0.$$

Solving this inequality, we obtain the existence of $t_{k+1} \geq t_k$ such that

$$\phi_{2^{k+1}}(t) \leq 2^k M_k^2, \qquad t \geq t_{k+1} \tag{2.13}$$

with

$$t_{k+1} - t_k \leq \frac{C}{2^k}. \tag{2.14}$$

Note that (2.12) and (2.13) imply

$$\phi_{2^{k+1}}(t) \leq M_{k+1}, \qquad t \geq t_{k+1}.$$

By the summability of the right side of (2.14) in $k$, the sequence $t_k$ with the indicated properties has been constructed. In particular,

$$\phi_{2^k} \leq M_k, \qquad t \geq T_0,$$

where $T_0 = \lim_k t_k < \infty$.

It remains to obtain a suitable upper bound for $M_k$. For this purpose, we construct a dominating sequence $R_1, R_2, R_3, \ldots$. Let

$$R_{k+1} = C_1 p_k^\mu R_k^2, \qquad k = 1, 2, \ldots \tag{2.15}$$

with a constant $C_1 \geq C_0$ to be determined and with $\mu = 2 + 2\lambda$. Also, set

$$R_1 = C_1 2^\mu M^2. \tag{2.16}$$

First, using induction, it is easy to check that (2.15) and (2.16) imply

$$R_k = (2^\mu C_1)^{2^k - 1} M^{2^k}, \qquad k = 1, 2, 3, \ldots. \tag{2.17}$$

Next, we claim that

$$M_k \leq R_k, \qquad k = 1, 2, \ldots. \tag{2.18}$$

It is clear that (2.18) holds for $k = 1$. Assuming that (2.18) holds for $k \in \mathbb{N}$, we get

$$\begin{aligned}
M_{k+1} &= C_0 \max\left\{ p_k M_k^2, p_k^{2(1+\lambda)p_k/(p_k+1)} M_k^{2p_k/(p_k+1)} M^{2p_k/(p_k+1)} \right\} \\
&\leq C_0 \max\left\{ p_k R_k^2, p_k^{2(1+\lambda)p_k/(p_k+1)} R_k^{2p_k/(p_k+1)} M^{2p_k/(p_k+1)} \right\} \\
&\leq C_1 p_k^\mu R_k^2 = R_{k+1}.
\end{aligned} \tag{2.19}$$

The second inequality in (2.19) is obtained by a direct verification. Since we have now established

$$M_k \leq (2^\mu C_1)^{2^k - 1} M^{2^k}, \qquad k = 1, 2, 3, \ldots,$$

by (2.17) and (2.18), we get

$$M_k^{1/2^k} \leq 2^\mu C_1 M, \qquad k = 1, 2, 3, \ldots,$$

and the first part of the lemma is established.

As for the last assertion, let $t_0 > 0$ be arbitrary. Applying the Gronwall lemma on (2.5), we get (2.6) for $t \geq t_0/2$, where $C$ depends on $\|u_0\|_{L^2}$ and $t_0$. By shifting time by $t_0/2$, we have (2.6) for $t \geq 0$. Similarly, we can choose $t_k = t_0/2^{k+1}$ for $k = 1, 2, \ldots$ and the constants then depend on $\|u_0\|_{L^2}$ and $t_0$. □

An important device in the proof of Theorem 2.1 is the modified vorticity

$$\zeta = \omega - R\rho, \tag{2.20}$$

introduced in [KW] where

$$R = \partial_1 \tilde{\Lambda}^{-2} = \partial_1 (I - \Delta)^{-1}$$

with $\tilde{\Lambda} = (I - \Delta)^{1/2}$. This, in turn, is a modification of the change of variable introduced in [JMWZ] (cf. also [SW, HKR]). The quantity $\zeta$ satisfies

$$\zeta_t - \Delta\zeta + u \cdot \nabla\zeta = [R, u \cdot \nabla]\rho - N\rho, \tag{2.21}$$

where

$$N = (\tilde{\Lambda}^{-2}\Delta - I)\partial_1 \tag{2.22}$$

is a smoothing operator of order $-1$ (cf. [KW]), i.e., the operator $\nabla N$ in the Calderón-Zygmund class. Using that $u$ is divergence-free, the first term on the right hand side of (2.21) may be rewritten as

$$[R, u \cdot \nabla]\rho = Ru_j \partial_j \rho - u_j \partial_j R\rho = \partial_j R(u_j \rho) - u_j \partial_j R\rho = [\partial_j R, u_j]\rho. \tag{2.23}$$

Also, for any multiplier operator $T$, we have

$$T\big([R, u \cdot \nabla]\rho\big) = [TR\partial_j, u_j]\rho - [T\partial_j, u_j]R\rho \tag{2.24}$$

(cf. [KW]). In both identities (2.23) and (2.24), which may be verified by a direct calculation, it is essential that $u$ is divergence-free.

*Proof of Theorem 2.1.* We assume

$$\|u_0\|_{H^2}, \|\rho_0\|_{H^1} \leq C. \tag{2.25}$$

By the Gagliardo-Nirenberg inequality

$$\|v\|_{L^p} \lesssim p^{1/2} \|v\|_{L^2}^{2/p} \|\nabla v\|_{L^2}^{1-2/p} + \|v\|_{L^2}$$

with $v = \rho_0$ and by (2.25), we get

$$\|\rho_0\|_{L^p} \lesssim p^{1/2}, \qquad p \in [2, \infty)$$

and thus

$$\|\rho(t)\|_{L^p} \lesssim p^{1/2}, \qquad t \geq 0, \qquad p \in [2, \infty). \tag{2.26}$$

Using (2.26) and applying Lemma 2.2 with $\lambda = 1/2$, there exists $t_1 \geq 0$ such that

$$\|\omega\|_{L^p} \leq C, \qquad t \geq t_1, \qquad p \in [2, \infty], \tag{2.27}$$

which by the triangle inequality implies

$$\|\zeta\|_{L^p} \lesssim 1, \qquad t \geq t_1, \qquad p \in [2, \infty]. \tag{2.28}$$

Since $C$ is allowed to depend on $\|u_0\|_{L^2}$, we may assume that $t_1 > 0$ is arbitrarily small.

In order to bound $\nabla \omega$, we consider evolution of the modified vorticity (2.20). Applying $\partial_k$ to (2.21), multiplying the resulting equation by $|\partial_k \zeta|^{2p-2} \partial_k \zeta$, integrating and summing in $k$ leads to

$$\begin{aligned}\frac{1}{2p}\frac{d}{dt}\sum_k \|\partial_k \zeta\|_{L^{2p}}^{2p} &- \sum_k \int (\Delta \partial_k \zeta)|\partial_k \zeta|^{2p-2}\partial_k \zeta\, dx\\ &= -\sum_k \int \partial_k(u_j \partial_j \zeta)|\partial_k \zeta|^{2p-2}\partial_k \zeta\, dx + \sum_k \int \partial_k([R, u\cdot\nabla]\rho)|\partial_k \zeta|^{2p-2}\partial_k \zeta\, dx\\ &\quad - \sum_k \int \partial_k N\rho|\partial_k \zeta|^{2p-2}\partial_k \zeta\, dx\\ &= J_1 + J_2 + J_3\end{aligned} \quad (2.29)$$

with no summation convention applied to the index $k$ in this proof. For $p \in \{2, 4, 8, \ldots\}$, denote

$$\psi_p = \sum_k \int |\partial_k \zeta|^p. \quad (2.30)$$

Note that the second term on the left hand side of (2.29) equals

$$\frac{2p-1}{p^2}\sum_k \int \partial_j(|\partial_k \zeta|^p)\partial_j(|\partial_k \zeta|^p) \geq \frac{D}{p},$$

where

$$D = \sum_k \int \partial_j(|\partial_k \zeta|^p)\partial_j(|\partial_k \zeta|^p) = \sum_k \|\nabla(|\partial_k \zeta|^p)\|_{L^2}^2.$$

Regarding $J_1$, we use the divergence-free condition on $u$ to write

$$\begin{aligned}J_1 &= -\sum_k \int \partial_k(u_j \partial_j \zeta)|\partial_k \zeta|^{2p-2}\partial_k \zeta = -\sum_k \int \partial_k u_j \partial_j \zeta |\partial_k \zeta|^{2p-2}\partial_k \zeta\\ &\lesssim \|\nabla u\|_{L^2}\|\nabla \zeta\|_{L^{4p}}\sum_k \||\partial_k \zeta|^{2p-1}\|_{L^{4p/(2p-1)}} \lesssim \|\nabla u\|_{L^2}\sum_k \|\nabla \zeta\|_{L^{4p}}^{2p}.\end{aligned}$$

Therefore,

$$J_1 \lesssim \|\omega\|_{L^2}\|\nabla \zeta\|_{L^{4p}}^{2p} \lesssim \sum_k \|\partial_k \zeta\|_{L^{4p}}^{2p}.$$

Using the Gagliardo-Nirenberg inequality, we have

$$\|\partial_k \zeta\|_{L^{4p}}^{2p} = \||\partial_k \zeta|^p\|_{L^4}^2 \lesssim \||\partial_k \zeta|^p\|_{L^2}\|\nabla(|\partial_k \zeta|^p)\|_{L^2} = \|\partial_k \zeta\|_{L^{2p}}^{p}\|\nabla(|\partial_k \zeta|^p)\|_{L^2}$$

for $k = 1, 2$, and thus

$$J_1 \leq \frac{D}{4p} + Cp\sum_k \|\partial_k \zeta\|_{L^{2p}}^{2p} \leq \frac{D}{4p} + Cp\psi_{2p}. \quad (2.31)$$

Next, for the second term $J_2$, we have

$$\begin{aligned}J_2 &= -(2p-1)\sum_k \int [R, u\cdot\nabla]\rho |\partial_k \zeta|^{2p-2}\partial_{kk} \zeta\, dx\\ &= -\frac{2p-1}{p}\sum_k \int [R, u\cdot\nabla]\rho |\partial_k \zeta|^{p-2}\partial_k \zeta \partial_k(|\partial_k \zeta|^p)\, dx\\ &\lesssim \|[R, u\cdot\nabla]\rho\|_{L^{2p}}\sum_k \|\nabla(|\partial_k \zeta|^p)\|_{L^2}\||\partial_k \zeta|^{p-1}\|_{L^{2p/(p-1)}}.\end{aligned}$$

The first factor is estimated as

$$\|[R, u \cdot \nabla]\rho\|_{L^{2p}} \le \|R(u_j \partial_j)\rho\|_{L^{2p}} + \|u_j \partial_j R\rho\|_{L^{2p}} = \|R\partial_j(u_j\rho)\|_{L^{2p}} + \|u_j \partial_j R\rho\|_{L^{2p}}$$
$$\lesssim p\|\rho u\|_{L^{2p}} + \|u\|_{L^\infty}\|(\nabla R)\rho\|_{L^{2p}} \lesssim p\|u\|_{L^\infty}\|\rho\|_{L^{2p}} \lesssim p^{3/2},$$

where we used (2.1) and (2.26) in the last inequality. Therefore, we obtain

$$\begin{aligned}
J_2 &\lesssim p^{3/2} \sum_k \|\nabla(|\partial_k\zeta|^p)\|_{L^2} \||\partial_k\zeta|^{p-1}\|_{L^{2p/(p-1)}} \lesssim p^{3/2} \sum_k \|\nabla(|\partial_k\zeta|^p)\|_{L^2} \|\partial_k\zeta\|_{L^{2p}}^{p-1} \\
&\le p^{3/2} D^{1/2} \sum_k \|\partial_k\zeta\|_{L^{2p}}^{p-1} \le \frac{D}{4p} + Cp^4 \sum_k \|\partial_k\zeta\|_{L^{2p}}^{2p-2} \le \frac{D}{4p} + Cp^4 \psi_{2p}^{(p-1)/p}.
\end{aligned} \quad (2.32)$$

For $J_3$, we use that the operator $N$, defined in (2.22), is a smoothing operator or order $-1$ (cf. [KW]). Thus

$$\begin{aligned}
J_3 &\lesssim \sum_k \|\partial_k N\rho\|_{L^{2p}} \||\partial_k\zeta|^{2p-1}\|_{L^{2p/(2p-1)}} \lesssim p \sum_k \|\rho\|_{L^{2p}} \||\partial_k\zeta|^{2p-1}\|_{L^{2p/(2p-1)}} \\
&\lesssim p^{3/2} \sum_k \|\partial_k\zeta\|_{L^{2p}}^{2p-1} \lesssim p^{3/2} \psi_{2p}^{(2p-1)/2p}.
\end{aligned} \quad (2.33)$$

By replacing the estimates (2.31), (2.32), and (2.33) in (2.29), we get

$$\frac{1}{p}\psi'_{2p} + \frac{1}{p}D \le Cp\psi_{2p} + Cp^4 \psi_{2p}^{(p-1)/p} + Cp^{3/2}\psi_{2p}^{(2p-1)/2p}, \qquad p \ge 2.$$

Using (2.9) with $v = |\partial_k\zeta|^p$, we obtain

$$\frac{1}{p}\psi'_{2p} + \frac{\psi_{2p}^2 - C\psi_p^4}{Cp\psi_p^2} \le Cp\psi_{2p} + Cp^4 \psi_{2p}^{(p-1)/p} + Cp^{3/2}\psi_{2p}^{(2p-1)/2p}, \qquad p \ge 2,$$

and thus, absorbing the last term on the right side and multiplying the resulting inequality by $p$,

$$\psi'_{2p} + \frac{\psi_{2p}^2}{C\psi_p^2} \le C\psi_p^2 + Cp^2 \psi_{2p} + Cp^5 \psi_{2p}^{(p-1)/p}, \qquad p \ge 2. \quad (2.34)$$

In order to start the induction, we also need an estimate for $\psi_2$. In this case, we have

$$D = \sum_k \int \partial_{jk}\zeta \partial_{jk}\zeta = \sum_k \|\nabla(\partial_k\zeta)\|_{L^2}^2 \gtrsim \frac{\|\nabla\zeta\|_{L^2}^4}{\|\zeta\|_{L^2}^2}.$$

Then the same derivation as above shows that

$$\psi'_2 + \frac{\|\nabla\zeta\|_{L^2}^4}{\|\zeta\|_{L^2}^2} \le C\psi_2 + C, \quad (2.35)$$

from where, using (2.28) with $p = 2$,

$$\psi'_2 + \psi_2^2 \le C\psi_2 + C.$$

Applying the Gronwall inequality, this implies that there exists $t_2 \ge t_1$ such that

$$\|\nabla\zeta\|_{L^2} \le C, \qquad t \ge t_2.$$

Going back to the inequality (2.34), fix $p \geq 2$, and note that if for any $t \geq 0$ we have

$$\psi_{2p} \geq C \max\{p^2 \psi_p^2, p^5 \psi_p^{2p/(p+1)}\},$$

for a sufficiently large constant $C$, half of the second term on the left hand side dominates the terms on the right hand side and thus

$$\psi'_{2p} + \frac{\psi_{2p}^2}{2\psi_p^2} \leq 0.$$

As in the proof of Lemma 2.2, this implies the existence of $t_3 \geq t_2$ such that

$$\|\nabla \zeta\|_{L^p} \leq C, \qquad t \geq t_3, \qquad p \in [2, \infty].$$

In particular, we get

$$\|\nabla \omega\|_{L^p} \lesssim p^{3/2}, \qquad t \geq t_3, \qquad p \in [2, \infty), \tag{2.36}$$

since $\|\nabla R\rho\|_{L^p} \lesssim p\|\rho\|_{L^p} \lesssim p^{3/2}$. The inequalities (2.27) and (2.36) then imply

$$\|\nabla u\|_{L^\infty} \leq C, \qquad t \geq t_3.$$

Since

$$\frac{d}{dt}\|\nabla \rho\|_{L^2}^2 \lesssim \|\nabla u\|_{L^\infty}\|\nabla \rho\|_{L^2}^2, \tag{2.37}$$

we get

$$\|\nabla \rho\|_{L^2} \leq Ce^{Ct}, \qquad t \geq 0,$$

and the assertion is proven. $\square$

*Remark* 2.3. It is not difficult to extend Theorem 2.1 to the case when we do not assume $\int_{\mathbb{T}^2} u_0 = 0$. In this case, we get $\left|\int_{\mathbb{T}^2} u\right| \lesssim t + 1$. Based on the energy inequality

$$\frac{1}{2}\frac{d}{dt}\|u\|_{L^2}^2 + \|\nabla u\|_{L^2}^2 \leq C\|u\|_{L^2}$$

we get $\|u(t)\|_{L^2} \lesssim t + 1$ for $t \geq 0$. Also, as in the proof above, we get $\|\omega\|_{L^p} \lesssim (t+1)^{1/2}$ for all $p \in [2, \infty]$ and thus also $\|\zeta\omega\|_{L^p} \lesssim (t+1)^{1/2}$ for all $t \geq t_1$ for some $t_1 \geq 0$. Again proceeding as above, we get $\|\nabla \zeta\|_{L^p} \lesssim (t+1)^{1/2}$ first for $p = 2$ and then for all $p \in [2, \infty]$ for $t$ sufficiently large.

# 3 The case $\mathbb{R}^2$

In this section, we consider the case of the whole space $\mathbb{R}^2$.

**Theorem 3.1.** *Assume that* $(u_0, \rho_0) \in H^2(\mathbb{R}^2) \times H^1(\mathbb{R}^2)$, *where* $\nabla \cdot u_0 = 0$. *Then we have*

$$\|u\|_{H^2} \leq C(t+1)^{1/2}, \qquad t \geq 0$$

and

$$\|\nabla\rho\|_{L^2} \leq Ce^{C(t+1)^{\beta+1}\log(t+2)}, \qquad t \geq 0$$

for a constant $C = C(\|u_0\|_{H^2}, \|\rho_0\|_{H^1})$, where

$$\beta = \prod_{j=1}^{\infty}\left(1 - \frac{1}{2^j}\right) = 0.28878\cdots.$$

Moreover,

$$\|\omega(t)\|_{L^p} \lesssim (t+1)^{1/p+\beta(1-2/p)}, \qquad t \geq t_0, \qquad p \in [2, \infty]$$

and

$$\|\nabla\omega(t)\|_{L^p} \lesssim p^{3/2} + (t+1)^{1/2}, \qquad t \geq t_0, \qquad p \in [2, \infty)$$

for some $t_0 \geq 0$.

*Remark* 3.2. The reason for a different bound than in Theorem 2.1 is a lack of the Poincaré inequality, which is available in other settings in this paper. If an additional damping term $\gamma u$, where $\gamma > 0$, is added to the left side of the equation (1.1), then the bounds are identical to those in Theorem 2.1, with constants depending on $\gamma$.

*Proof of Theorem 3.1.* The energy inequality

$$\frac{1}{2}\frac{d}{dt}\|u\|_{L^2}^2 + \|\nabla u\|_{L^2}^2 \leq C\|u\|_{L^2}$$

implies

$$\|u(t)\|_{L^2} \lesssim t+1, \qquad t \geq 0.$$

Similarly, the $L^2$ inequality for the vorticity reads

$$\frac{1}{2}\frac{d}{dt}\|\omega\|_{L^2}^2 + \|\nabla\omega\|_{L^2}^2 \leq C,$$

which implies

$$\|\nabla u(t)\|_{L^2} = \|\omega\|_{L^2} \lesssim (t+1)^{1/2}, \qquad t \geq 0. \tag{3.1}$$

Next, we consider the upper bounds for $\|\omega\|_{L^p}$ and $\|\nabla\omega\|_{L^p}$ for $p \geq 2$. Denote

$$\phi_p = \|\omega\|_{L^p}^p$$

and fix $p \geq 2$. From the vorticity equation

$$\omega_t - \Delta\omega + u \cdot \nabla\omega = \partial_1\rho$$

we obtain, as in (2.10), the inequality

$$\phi'_{2p} + \frac{\phi_{2p}^2}{C\phi_p^2} \leq Cp^3\phi_{2p}^{(p-1)/p}.$$

As in the proof of Lemma 2.2, we conclude by induction that

$$\|\omega\|_{L^p} \lesssim (t+1)^{\beta_p}, \qquad t \geq t_0,$$

for $p = 2, 4, \ldots$, where

$$\beta_{2^k} = \prod_{j=1}^{k}\left(1 - \frac{1}{2^j}\right).$$

Therefore,

$$\|\omega\|_{L^\infty} \lesssim (t+1)^{\beta}, \qquad t \geq t_0.$$

Combined with (3.1), we get

$$\|\omega\|_{L^p} \lesssim (t+1)^{1/p+\beta(1-2/p)}, \qquad t \geq t_0, \qquad p \in [2, \infty],$$

from where also

$$\|\nabla u\|_{L^p} \lesssim p(t+1)^{1/p+\beta(1-2/p)}, \qquad t \geq t_0, \qquad p \in [2, \infty).$$

In order to obtain an estimate on the growth of $\nabla \omega$, we consider the generalized vorticity (2.20), which satisfies (2.21). As in the periodic case, we set (2.30), i.e.,

$$\psi_p = \sum_k \int |\partial_k \zeta|^p, \qquad p \geq 2$$

and obtain

$$\psi_2' + \frac{\psi_2^2}{\|\zeta\|_{L^2}^2} \leq C\psi_2 + C \tag{3.2}$$

(cf. (2.35) above) and

$$\psi_{2p}' + \frac{\psi_{2p}^2}{C\psi_p^2} \leq Cp^2 \psi_{2p} + Cp^5 \psi_{2p}^{(p-1)/p}, \qquad p \geq 2 \tag{3.3}$$

(cf. (2.34) above). The inequality (3.2) and

$$\|\zeta\|_{L^2} \lesssim (t+1)^{1/2}, \qquad t \geq 0$$

imply

$$\psi_2(t) \lesssim t+1, \qquad t \geq 0.$$

Continuing by induction, we obtain from (3.3)

$$\psi_p(t) \lesssim p^\mu (t+1)^{p/2}, \qquad t \geq 0, \qquad p = 2, 2^2, 2^3, \ldots$$

with a certain $\mu > 0$. These inequalities then lead to

$$\|\nabla \zeta\|_{L^p} \lesssim (t+1)^{1/2}, \qquad t \geq 0, \qquad p \in [2, \infty].$$

From here, we obtain $\|\nabla \omega\|_{L^p} \leq \|\nabla \zeta\|_{L^p} + \|\nabla R\rho\|_{L^p} \lesssim (t+1)^{1/2} + p^{3/2}$, and thus $\|D^2 u\|_{L^p} \lesssim p^{5/2}(t+1)^{1/2}$. Therefore,

$$\|\nabla u\|_{L^\infty} \leq C\|\nabla u\|_{L^p}^{1-2/p}\|D^2 u\|_{L^p}^{2/p} \lesssim Cp(t+1)^{3/p-2/p^2}.$$

Choosing a proper value for $p$, we get

$$\|\nabla u\|_{L^\infty} \lesssim (t+1)^{\beta_\infty} \log(t+1), \qquad t \geq t_0$$

which then implies

$$\|\nabla \rho\|_{L^2} \lesssim \exp\left((t+1)^{\beta_\infty+1} \log(t+1)\right), \qquad t \geq 0,$$

and the theorem is proven. $\square$

## 4 Bounds with the Lions boundary condition

In this section, we consider the Boussinesq system on a bounded smooth domain $\Omega \subseteq \mathbb{R}^2$, with the Lions boundary conditions

$$u \cdot n = \omega = 0 \quad \text{on } \partial\Omega,$$

where $n$ denotes the outward unit normal. We use the standard notation corresponding to the Navier-Stokes system [CF, T1, R, HKZ1]. In particular, denote

$$H = \{u \in L^2(\Omega) : \nabla \cdot u = 0, u \cdot n = 0 \text{ on } \partial\Omega\},$$

where $n$ stands for the outward unit normal vector with respect to the domain $\Omega$, which is assumed to be smooth and bounded. Let also

$$V = \{u \in H^1(\Omega) : \nabla \cdot u = 0, u \cdot n = 0 \text{ on } \partial\Omega\}.$$

The Stokes operator $A \colon D(A) \to H$, with the domain $D(A) = H^2(\Omega) \cap V$, is defined by $A = -\mathbb{P}\Delta$, where $\mathbb{P}$ is the Leray projector in $L^2(\Omega)$ on the space $H$.

**Theorem 4.1.** *Assume that $(u_0, \rho_0) \in D(A) \times H^1(\Omega)$. Then we have*

$$\|u\|_{H^2} \leq C, \qquad t \geq 0 \tag{4.1}$$

*and*

$$\|\nabla \rho\|_{L^2} \leq C e^{Ct}, \qquad t \geq 0 \tag{4.2}$$

*for a constant $C = C(\|u_0\|_{D(A)}, \|\rho_0\|_{H^1})$. In addition, we have*

$$\|\omega(t)\|_{L^p} \leq C, \qquad t \geq t_0, \qquad p \in [2, \infty],$$

*where $t_0 \geq 0$ depends on $\|u_0\|_{L^2}$ and $\|\rho_0\|_{L^2}$.*

The global persistence for the Boussinesq system with the Lions boundary conditions was recently addressed by Doering et al in [DWZZ]. The authors moreover proved that $\|u\|_{H^1} \to 0$ as $t \to \infty$. It is not clear whether the same holds for other boundary conditions considered in the present paper. Namely, the important ingredients in [DWZZ] are that $\theta = ay + b$ belongs to the state space and that the vorticity $\omega$ vanishes on the boundary.

From here on, the constant $C$ is allowed to depend on $\|u_0\|_{D(A)}$ and $\|\rho_0\|_{H^1}$. The proof of the assertion (4.1) is the same as in [J], which considered the Dirichlet boundary condition. From [J], we also recall the inequality
$$\int_{t_1}^{t_2} \|A^{3/2}u(s)\|_{L^2}^2 \, ds \leq C(t_2 - t_1 + 1), \qquad 0 \leq t_1 \leq t_2$$
(cf. [J, p. 115]).

*Proof of Theorem 4.1.* Note that the proof of Lemma 2.2 applies here verbatim, and thus we obtain
$$\|\omega(\cdot, t)\|_{L^p} \leq C, \qquad t \geq t_0, \qquad 2 \leq p \leq \infty. \tag{4.3}$$
Since $t_0 > 0$ may be chosen arbitrarily small (cf. Lemma 2.2) and by the local existence, we may simply assume that (4.3) holds for all $t \geq 0$.

Now, note that the argument starting in (2.29) does not apply in this setting due to arising boundary terms. Thus we use an alternative argument, described next. Fix $t_0 > 0$. Let $\theta \colon \mathbb{R} \to [0, \infty)$ be a smooth non-decreasing function such that $\theta \equiv 0$ on $[0, t_0/2]$ and $\theta \equiv 1$ on $[t_0, \infty]$. Then we have
$$\partial_t(\theta(t)\omega) - \Delta(\theta(t)\omega) = \theta'(t)\omega - \partial_j(\theta(t)u_j\omega) + \partial_1(\theta(t)\rho)$$
$$= \partial_1(\theta'(t)u_2) - \partial_2(\theta'(t)u_1) - \partial_j(\theta(t)u_j\omega) + \partial_1(\theta(t)\rho).$$
Using the parabolic regularity with the right side in divergence form we get, for all $t \geq 0$,
$$\begin{aligned}
\left(\int_0^t \|\theta(s)\nabla\omega(s)\|_{L^p}^p \, ds\right)^{1/p} &\\
&\hspace{-3cm}\leq Cp\left(\int_0^t \|\theta'(s)u\|_{L^p}^p \, ds\right)^{1/p} + Cp\left(\int_0^t \|\theta(s)\omega(s)u(s)\|_{L^p}^p \, ds\right)^{1/p} + Cp\left(\int_0^t \|\theta(s)\rho(s)\|_{L^p}^p \, ds\right)^{1/p} \\
&\hspace{-3cm}\leq Cp\left(\int_0^t \|\theta'(s)u\|_{L^p}^p \, ds\right)^{1/p} + Cp\left(\int_0^t \|\theta(s)\omega(s)\|_{L^p}^p \, ds\right)^{1/p} + Cp\left(\int_0^t \|\theta(s)\rho(s)\|_{L^p}^p \, ds\right)^{1/p} \\
&\hspace{-3cm}\leq Cp\left(\int_0^t \|u\|_{L^p}^p \, ds\right)^{1/p} + Cp^{3/2}t^{1/p},
\end{aligned} \tag{4.4}$$
where $2 \leq p < \infty$ by $\|u\|_{L^\infty} \lesssim \|u\|_{H^2} \lesssim 1$. Therefore, using $\|u\|_{L^p} \lesssim 1$,
$$\left(\int_0^t \|\theta(s)\nabla\omega(s)\|_{L^p}^p \, ds\right)^{1/p} \leq Cp^{3/2}t^{1/p}.$$
Now, for every $p \in [2, \infty)$, we have
$$\|\nabla u\|_{L^\infty} \leq C\|\nabla u\|_{L^p}^{1-2/p}\|D^2 u\|_{L^p}^{2/p} + C\|\nabla u\|_{L^p} \leq Cp\|\omega\|_{L^p}^{1-2/p}\|\nabla\omega\|_{L^p}^{2/p} + Cp\|\omega\|_{L^p}.$$
In particular,
$$\begin{aligned}
\int_{t_0}^t \|\nabla u\|_{L^\infty} \, ds &\lesssim \int_{t_0}^t \|\omega\|_{L^4}^{1/2}\|\nabla\omega\|_{L^4}^{1/2} \, ds + \int_{t_0}^t \|\omega\|_{L^4} \, ds \\
&\lesssim \left(\int_{t_0}^t \|\omega\|_{L^4}^{4/7} \, ds\right)^{7/8}\left(\int_{t_0}^t \|\nabla\omega\|_{L^4}^4 \, ds\right)^{1/8} + \int_{t_0}^t \|\omega\|_{L^4} \, ds \\
&\lesssim t^{7/8}t^{1/8} + t \lesssim t,
\end{aligned} \tag{4.5}$$

where we used (4.4) with $p = 4$ in the last inequality. Integrating (2.37), which also holds in this setting, and applying (4.5) then gives the inequality (4.2). □

## 5 Bounds with the Dirichlet boundary condition

Finally, we address the long time behavior of the Boussinesq system with the classical Dirichlet (non-slip) boundary condition

$$u = 0 \quad \text{on } \partial\Omega,$$

where $\Omega$ is a bounded smooth domain. Recall the standard notation $H = \{u \in L^2(\Omega) : \nabla \cdot u = 0, u \cdot n = 0 \text{ on } \partial\Omega\}$, where $n$ denotes the outward unit normal vector with respect to the domain $\Omega$, and $V = H_0^1(\Omega) \cap H$. The Stokes operator is then defined as in the previous section, i.e.,

$$A = -\mathbb{P}\Delta,$$

with the domain $D(A) = H^2(\Omega) \cap V$, where $\mathbb{P}$ is the Leray projector in $L^2(\Omega)$ on the space $H$.

**Theorem 5.1.** *Assume that $(u_0, \rho_0) \in D(A) \times H^1(\Omega)$. Then we have*

$$\|v\|_{H^2} \leq C, \qquad t \geq 0 \tag{5.1}$$

*and*

$$\|\nabla \rho\|_{L^2} \leq Ce^{Ct}, \qquad t \geq 0 \tag{5.2}$$

*for a constant $C = C(\|u_0\|_{D(A)}, \|\rho_0\|_{H^1})$.*

*Proof of Theorem 5.1.* With $\theta = \theta(t)$ a smooth cut-off function as in the previous section, we have

$$\partial_t(\theta u) - \Delta(\theta u) + u \cdot \nabla(\theta u) + \nabla(\theta p) = \theta' u + \theta \rho e_2.$$

Using the $W^{2,4}$ regularity estimate due to Sohr and Von Wahl [SvW], we get

$$\left(\int_0^t \|\theta D^2 u\|_{L^4}^4\right)^{1/4} \lesssim \left(\int_0^t \|u \cdot \nabla(\theta u)\|_{L^4}^4 \, ds\right)^{1/4} + \left(\int_0^t \|\theta' u\|_{L^4}^4 \, ds\right)^{1/4} + \left(\int_0^t \|\rho\|_{L^4}^4 \, ds\right)^{1/4}$$

$$\lesssim \left(\int_0^t \|u\|_{L^8}^4 \|\nabla u\|_{L^8}^4 \, ds\right)^{1/4} + \left(\int_0^t \|u\|_{L^4}^4 \, ds\right)^{1/2} + t^{1/4}$$

$$\lesssim \left(\int_0^t \|u\|_{L^2} \|\nabla u\|_{L^2}^4 \|D^2 u\|_{L^2}^3 \, ds\right)^{1/4} + \left(\int_0^t \|u\|_{L^2}^2 \|\nabla u\|_{L^2}^2 \, ds\right)^{1/4} + t^{1/4}$$

$$\lesssim t^{1/4}$$

whence

$$\int_{t_0}^t \|D^2 u\|_{L^4}^4 \lesssim t. \tag{5.3}$$

Also, by (5.1), we obtain

$$\|\omega(t)\|_{L^p} \leq C(p), \qquad t \geq t_0, \qquad p \in [2, \infty). \tag{5.4}$$

As in the previous section, the inequalities (5.3) and (5.4) with $p = 4$ imply

$$\int_{t_0}^t \|\nabla u\|_{L^\infty} \, ds \lesssim t, \qquad t \geq t_0,$$

and (5.2) follows from (2.37). □

# Acknowledgments

The authors were supported in part by the NSF grants DMS-1615239 and DMS-1907992.